# SOME COMPANIONS OF GRÜSS INEQUALITY IN INNER PRODUCT SPACES

S.S. DRAGOMIR

ABSTRACT. Some companions of Grüss inequality in inner product spaces and applications for integrals are given.

## 1. INTRODUCTION

The following inequality of Grüss type in real or complex linear spaces is known (see [1]).

**Theorem 1.** *Let $(H; \langle \cdot, \cdot \rangle)$ be an inner product space over $\mathbb{K}$ ($\mathbb{K} = \mathbb{C}, \mathbb{R}$) and $e \in H$, $\|e\| = 1$. If $\phi, \gamma, \Phi, \Gamma$ are real or complex numbers and $x, y$ are vectors in $H$ such that the condition*

$$(1.1) \qquad \operatorname{Re} \langle \Phi e - x, x - \phi e \rangle \geq 0 \quad \text{and} \quad \operatorname{Re} \langle \Gamma e - y, y - \gamma e \rangle \geq 0$$

*or, equivalently (see [2]),*

$$(1.2) \qquad \left\| x - \frac{\phi + \Phi}{2} e \right\| \leq \frac{1}{2} |\Phi - \phi| \quad \text{and} \quad \left\| y - \frac{\gamma + \Gamma}{2} e \right\| \leq \frac{1}{2} |\Gamma - \gamma|$$

*holds, then we have the inequality*

$$(1.3) \qquad |\langle x, y \rangle - \langle x, e \rangle \langle e, y \rangle| \leq \frac{1}{4} |\Phi - \phi| |\Gamma - \gamma|.$$

*The constant $\frac{1}{4}$ is best possible in the sense that it cannot be replaced by a smaller constant.*

The following particular instances for integrals and means are useful in applications.

**Corollary 1.** *Let $f, g : [a, b] \to \mathbb{K}$ ($\mathbb{K} = \mathbb{C}, \mathbb{R}$) be Lebesgue measurable and such that there exists the constants $\phi, \gamma, \Phi, \Gamma \in \mathbb{K}$ with the property*

$$(1.4) \qquad \operatorname{Re}\left[(\Phi - f(x))\left(\overline{f(x)} - \overline{\phi}\right)\right] \geq 0, \quad \operatorname{Re}\left[(\Gamma - g(x))\left(\overline{g(x)} - \overline{\gamma}\right)\right] \geq 0$$

*for a.e. $x \in [a, b]$, or, equivalently*

$$(1.5) \qquad \left| f(x) - \frac{\phi + \Phi}{2} \right| \leq \frac{1}{2} |\Phi - \phi| \quad \text{and} \quad \left| g(x) - \frac{\gamma + \Gamma}{2} \right| \leq \frac{1}{2} |\Gamma - \gamma|$$

*for a.e. $x \in [a, b]$.*

*Then we have the inequality*

$$\left| \frac{1}{b-a} \int_a^b f(x) \overline{g(x)} dx - \frac{1}{b-a} \int_a^b f(x) dx \cdot \frac{1}{b-a} \int_a^b \overline{g(x)} dx \right| \tag{1.6}$$

$$\leq \frac{1}{4} |\Phi - \phi| |\Gamma - \gamma|.$$

*The constant $\frac{1}{4}$ is best possible.*

The discrete case is incorporated in

**Corollary 2.** *Let $\mathbf{x}, \mathbf{y} \in \mathbb{K}^n$ and $\phi, \gamma, \Phi, \Gamma \in \mathbb{K}$ be such that*

$$\operatorname{Re}\left[(\Phi - x_i)\left(\overline{x_i} - \overline{\phi}\right)\right] \geq 0 \quad \text{and} \quad \operatorname{Re}\left[(\Gamma - y_i)\left(\overline{y_i} - \overline{\gamma}\right)\right] \geq 0, \tag{1.7}$$

*for each $i \in \{1, \ldots, n\}$, or, equivalently,*

$$\left| x_i - \frac{\phi + \Phi}{2} \right| \leq \frac{1}{2} |\Phi - \phi| \quad \text{and} \quad \left| y_i - \frac{\gamma + \Gamma}{2} \right| \leq \frac{1}{2} |\Gamma - \gamma|, \tag{1.8}$$

*for each $i \in \{1, \ldots, n\}$.*

*Then we have the inequality*

$$\left| \frac{1}{n} \sum_{i=1}^n x_i \overline{y_i} - \frac{1}{n} \sum_{i=1}^n x_i \cdot \frac{1}{n} \sum_{i=1}^n \overline{y_i} \right| \leq \frac{1}{4} |\Phi - \phi| |\Gamma - \gamma|. \tag{1.9}$$

*The constant $\frac{1}{4}$ is best possible in (1.9).*

For some recent results related to Grüss type inequalities in inner product spaces, see [2]. More applications of Theorem 1 for integral and discrete inequalities may be found in [3].

It is the main aim of this paper to provide other inequalities of Grüss type in the general seting of inner product spaces over the real or complex number field $\mathbb{K}$. Applications for Lebsegue integrals are pointed out as well.

## 2. A Grüss Type Inequality

The following Grüss type inequality in inner product spaces holds.

**Theorem 2.** *Let $x, y, e \in H$ with $\|e\| = 1$, and the scalars $a, A, b, B \in \mathbb{K}$ ($\mathbb{K} = \mathbb{C}, \mathbb{R}$) such that $\operatorname{Re}(\bar{a}A) > 0$ and $\operatorname{Re}(\bar{b}B) > 0$. If*

$$\operatorname{Re}\langle Ae - x, x - ae \rangle \geq 0 \quad \text{and} \quad \operatorname{Re}\langle Be - y, y - be \rangle \geq 0 \tag{2.1}$$

*or, equivalently (see [2]),*

$$\left\| x - \frac{a+A}{2} e \right\| \leq \frac{1}{2} |A - a| \quad \text{and} \quad \left\| y - \frac{b+B}{2} e \right\| \leq \frac{1}{2} |B - b|, \tag{2.2}$$

*then we have the inequality*

$$|\langle x, y \rangle - \langle x, e \rangle \langle e, y \rangle| \leq \frac{1}{4} M(a, A) M(b, B) |\langle x, e \rangle \langle e, y \rangle|, \tag{2.3}$$

*where $M(\cdot, \cdot)$ is defined by*

$$M(a, A) := \left[ \frac{(|A| - |a|)^2 + 4\left[|A\bar{a}| - \operatorname{Re}(A\bar{a})\right]}{\operatorname{Re}(\bar{a}A)} \right]^{\frac{1}{2}}. \tag{2.4}$$

The constant $\frac{1}{4}$ is best possible in the sense that it cannot be replaced by a smaller constant.

*Proof.* Apply Schwartz's inequality in $(H; \langle \cdot, \cdot \rangle)$ for the vectors $x - \langle x, e \rangle e$ and $y - \langle y, e \rangle e$, to get (see also [1])

$$(2.5) \qquad |\langle x, y \rangle - \langle x, e \rangle \langle e, y \rangle|^2 \leq \left( \|x\|^2 - |\langle x, e \rangle|^2 \right) \left( \|y\|^2 - |\langle y, e \rangle|^2 \right).$$

Now, assume that $u, v \in H$, and $c, C \in \mathbb{K}$ with $\operatorname{Re}(\bar{c}C) > 0$ and $\operatorname{Re} \langle Cv - u, u - cv \rangle \geq 0$. This last inequality is equivalent to

$$(2.6) \qquad \|u\|^2 + \operatorname{Re}(\bar{c}C) \|v\|^2 \leq \operatorname{Re}\left[ C\overline{\langle u, v \rangle} + \bar{c} \langle u, v \rangle \right].$$

Dividing this inequality by $[\operatorname{Re}(C\bar{c})]^{\frac{1}{2}} > 0$, we deduce

$$(2.7) \qquad \frac{1}{[\operatorname{Re}(\bar{c}C)]^{\frac{1}{2}}} \|u\|^2 + [\operatorname{Re}(\bar{c}C)]^{\frac{1}{2}} \|v\|^2 \leq \frac{\operatorname{Re}\left[ C\overline{\langle u, v \rangle} + \bar{c} \langle u, v \rangle \right]}{[\operatorname{Re}(\bar{c}C)]^{\frac{1}{2}}}.$$

On the other hand, by the elementary inequality

$$\alpha p^2 + \frac{1}{\alpha} q^2 \geq 2pq, \quad \alpha > 0, \ p, q \geq 0,$$

we deduce

$$(2.8) \qquad 2 \|u\| \|v\| \leq \frac{1}{[\operatorname{Re}(\bar{c}C)]^{\frac{1}{2}}} \|u\|^2 + [\operatorname{Re}(\bar{c}C)]^{\frac{1}{2}} \|v\|^2.$$

Making use of (2.7) and (2.8) and the fact that for any $z \in \mathbb{C}$, $\operatorname{Re}(z) \leq |z|$, we get

$$\|u\| \|v\| \leq \frac{\operatorname{Re}\left[ C\overline{\langle u, v \rangle} + \bar{c} \langle u, v \rangle \right]}{2 [\operatorname{Re}(\bar{c}C)]^{\frac{1}{2}}} \leq \frac{|c| + |C|}{2 [\operatorname{Re}(\bar{c}C)]^{\frac{1}{2}}} |\langle u, v \rangle|.$$

Consequently

$$(2.9) \qquad \|u\|^2 \|v\|^2 - |\langle u, v \rangle|^2 \leq \left[ \frac{(|c| + |C|)^2}{4 [\operatorname{Re}(\bar{c}C)]} - 1 \right] |\langle u, v \rangle|^2$$

$$= \frac{1}{4} \frac{(|c| - |C|)^2 + 4 [|\bar{c}C| - \operatorname{Re}(\bar{c}C)]}{\operatorname{Re}(\bar{c}C)} |\langle u, v \rangle|^2$$

$$= \frac{1}{4} M^2(c, C) |\langle u, v \rangle|^2.$$

Now, if we write (2.9) for the choices $u = x$, $v = e$ and $u = y$, $v = e$ respectively and use (2.5), we deduce the desired result (2.2). The sharpness of the constant will be proved in the case where $H$ is a real inner product space. ∎

The following corollary which provides a simpler Grüss type inequality for real constants (and in particular, for real inner product spaces) holds.

**Corollary 3.** *With the assumptions of Theorem 2 and if $a, b, A, B \in \mathbb{R}$ are such that $A > a > 0$, $B > b > 0$ and*

$$(2.10) \qquad \left\| x - \frac{a + A}{2} e \right\| \leq \frac{1}{2}(A - a) \quad \text{and} \quad \left\| y - \frac{b + B}{2} e \right\| \leq \frac{1}{2}(B - b),$$

*then we have the inequality*

$$(2.11) \quad |\langle x, y \rangle - \langle x, e \rangle \langle e, y \rangle| \leq \frac{1}{4} \cdot \frac{(A-a)(B-b)}{\sqrt{abAB}} |\langle x, e \rangle \langle e, y \rangle|.$$

*The constant $\frac{1}{4}$ is best possible.*

*Proof.* The prove the sharpness of the constant $\frac{1}{4}$ assume that the inequality (2.11) holds in real inner product spaces with $x = y$ and for a constant $k > 0$, i.e.,

$$(2.12) \quad \|x\|^2 - |\langle x, e \rangle|^2 \leq k \cdot \frac{(A-a)^2}{aA} |\langle x, e \rangle|^2 \quad (A > a > 0),$$

provided $\left\| x - \frac{a+A}{2} e \right\| \leq \frac{1}{2}(A-a)$, or equivalently, $\langle Ae - x, x - ae \rangle \geq 0$.

We choose $H = \mathbb{R}^2$, $x = (x_1, x_2) \in \mathbb{R}^2$, $e = \left( \frac{1}{\sqrt{2}}, \frac{1}{\sqrt{2}} \right)$. Then we have

$$\|x\|^2 - |\langle x, e \rangle|^2 = x_1^2 + x_2^2 - \frac{(x_1 + x_2)^2}{2} = \frac{(x_1 - x_2)^2}{2},$$

$$|\langle x, e \rangle|^2 = \frac{(x_1 + x_2)^2}{2},$$

and by (2.12) we get

$$(2.13) \quad \frac{(x_1 - x_2)^2}{2} \leq k \cdot \frac{(A-a)^2}{aA} \cdot \frac{(x_1 + x_2)^2}{2}.$$

Now, if we let $x_1 = \frac{a}{\sqrt{2}}$, $x_2 = \frac{A}{\sqrt{2}}$ $(A > a > 0)$, then obviously

$$\langle Ae - x, x - ae \rangle = \sum_{i=1}^{2} \left( \frac{A}{\sqrt{2}} - x_i \right) \left( x_i - \frac{a}{\sqrt{2}} \right) = 0,$$

which shows that the condition (2.10) is fulfilled, and by (2.13) we get

$$\frac{(A-a)^2}{4} \leq k \cdot \frac{(A-a)^2}{aA} \cdot \frac{(a+A)^2}{4}$$

for any $A > a > 0$. This implies

$$(2.14) \quad (A+a)^2 k \geq aA$$

for any $A > a > 0$.

Finally, let $a = 1 - q$, $A = 1 + q$, $q \in (0, 1)$. Then from (2.14) we get $4k \geq 1 - q^2$ for any $q \in (0, 1)$ which produces $k \geq \frac{1}{4}$. ∎

**Remark 1.** *If $\langle x, e \rangle, \langle y, e \rangle$ are assumed not to be zero, then the inequality (2.3) is equivalent to*

$$(2.15) \quad \left| \frac{\langle x, y \rangle}{\langle x, e \rangle \langle e, y \rangle} - 1 \right| \leq \frac{1}{4} M(a, A) M(b, B),$$

*while the inequality (2.11) is equivalent to*

$$(2.16) \quad \left| \frac{\langle x, y \rangle}{\langle x, e \rangle \langle e, y \rangle} - 1 \right| \leq \frac{1}{4} \cdot \frac{(A-a)(B-b)}{\sqrt{abAB}}.$$

*The constant $\frac{1}{4}$ is best possible in both inequalities.*

## 3. Some Related Results

The following result holds.

**Theorem 3.** *Let $(H; \langle \cdot, \cdot \rangle)$ be an inner product space over $\mathbb{K}$ ($\mathbb{K} = \mathbb{C}, \mathbb{R}$). If $\gamma, \Gamma \in \mathbb{K}$, $e, x, y \in H$ with $\|e\| = 1$ and $\lambda \in (0, 1)$ are such that*

$$(3.1) \quad \operatorname{Re} \langle \Gamma e - (\lambda x + (1 - \lambda) y), (\lambda x + (1 - \lambda) y) - \gamma e \rangle \geq 0,$$

*or, equivalently,*

$$(3.2) \quad \left\| \lambda x + (1 - \lambda) y - \frac{\gamma + \Gamma}{2} e \right\| \leq \frac{1}{2} |\Gamma - \gamma|,$$

*then we have the inequality*

$$(3.3) \quad \operatorname{Re} [\langle x, y \rangle - \langle x, e \rangle \langle e, y \rangle] \leq \frac{1}{16} \cdot \frac{1}{\lambda (1 - \lambda)} |\Gamma - \gamma|^2.$$

*The constant $\frac{1}{16}$ is the best possible constant in (3.3) in the sense that it cannot be replaced by a smaller one.*

*Proof.* We know that for any $z, u \in H$ one has

$$\operatorname{Re} \langle z, u \rangle \leq \frac{1}{4} \|z + u\|^2.$$

Then for any $a, b \in H$ and $\lambda \in (0, 1)$ one has

$$(3.4) \quad \operatorname{Re} \langle a, b \rangle \leq \frac{1}{4\lambda (1 - \lambda)} \|\lambda a + (1 - \lambda) b\|^2.$$

Since

$$\langle x, y \rangle - \langle x, e \rangle \langle e, y \rangle = \langle x - \langle x, e \rangle e, y - \langle y, e \rangle e \rangle \quad (\text{as } \|e\| = 1),$$

using (3.4), we have

$$(3.5) \quad \operatorname{Re} [\langle x, y \rangle - \langle x, e \rangle \langle e, y \rangle]$$
$$= \operatorname{Re} [\langle x - \langle x, e \rangle e, y - \langle y, e \rangle e \rangle]$$
$$\leq \frac{1}{4\lambda (1 - \lambda)} \|\lambda (x - \langle x, e \rangle e) + (1 - \lambda) (y - \langle y, e \rangle e)\|^2$$
$$= \frac{1}{4\lambda (1 - \lambda)} \|\lambda x + (1 - \lambda) y - \langle \lambda x + (1 - \lambda) y, e \rangle e\|^2.$$

Since, for $m, e \in H$ with $\|e\| = 1$, one has the equality

$$(3.6) \quad \|m - \langle m, e \rangle e\|^2 = \|m\|^2 - |\langle m, e \rangle|^2,$$

then by (3.5) we deduce the inequality

$$(3.7) \quad \operatorname{Re} [\langle x, y \rangle - \langle x, e \rangle \langle e, y \rangle]$$
$$\leq \frac{1}{4\lambda (1 - \lambda)} \left[ \|\lambda x + (1 - \lambda) y\|^2 - |\langle \lambda x + (1 - \lambda) y, e \rangle|^2 \right].$$

Now, if we apply Grüss' inequality

$$0 \leq \|a\|^2 - |\langle a, e \rangle|^2 \leq \frac{1}{4} |\Gamma - \gamma|^2$$

provided

$$\operatorname{Re} \langle \Gamma e - a, a - \gamma e \rangle \geq 0,$$

for $a = \lambda x + (1 - \lambda) y$, we have

(3.8) $$\|\lambda x + (1 - \lambda) y\|^2 - |\langle \lambda x + (1 - \lambda) y, e \rangle|^2 \leq \frac{1}{4} |\Gamma - \gamma|^2.$$

Utilising (3.7) and (3.8) we deduce the desired inequality (3.3). To prove the sharpness of the constant $\frac{1}{16}$, assume that (3.3) holds with a constant $C > 0$, provided (3.1) is valid, i.e.,

(3.9) $$\operatorname{Re}\left[\langle x, y \rangle - \langle x, e \rangle \langle e, y \rangle\right] \leq C \cdot \frac{1}{\lambda (1 - \lambda)} |\Gamma - \gamma|^2.$$

If in (3.9) we choose $x = y$, provided (3.1) holds with $x = y$ and $\lambda \in (0, 1)$, then

(3.10) $$\|x\|^2 - |\langle x, e \rangle|^2 \leq C \cdot \frac{1}{\lambda (1 - \lambda)} |\Gamma - \gamma|^2,$$

provided
$$\operatorname{Re} \langle \Gamma e - x, x - \gamma e \rangle \geq 0.$$

Since we know, in Grüss' inequality, the constant $\frac{1}{4}$ is best possible, then by (3.10), one has
$$\frac{1}{4} \leq \frac{C}{\lambda (1 - \lambda)} \quad \text{for} \quad \lambda \in (0, 1),$$

giving, for $\lambda = \frac{1}{2}$, $C \geq \frac{1}{16}$.

The theorem is completely proved. ∎

The following corollary is a natural consequence of the above result.

**Corollary 4.** *Assume that $\gamma, \Gamma, e, x, y$ and $\lambda$ are as in Theorem 3. If*

(3.11) $$\operatorname{Re} \langle \Gamma e - (\lambda x \pm (1 - \lambda) y), (\lambda x \pm (1 - \lambda) y) - \gamma e \rangle \geq 0,$$

*or, equivalently,*

(3.12) $$\left\|\lambda x \pm (1 - \lambda) y - \frac{\gamma + \Gamma}{2} e\right\| \leq \frac{1}{2} |\Gamma - \gamma|^2,$$

*then we have the inequality*

(3.13) $$|\operatorname{Re}\left[\langle x, y \rangle - \langle x, e \rangle \langle e, y \rangle\right]| \leq \frac{1}{16} \cdot \frac{1}{\lambda (1 - \lambda)} |\Gamma - \gamma|^2.$$

*The constant $\frac{1}{16}$ is best possible in (3.13).*

*Proof.* Using Theorem 3 for $(-y)$ instead of $y$, we have that
$$\operatorname{Re} \langle \Gamma e - (\lambda x - (1 - \lambda) y), (\lambda x - (1 - \lambda) y) - \gamma e \rangle \geq 0,$$

which implies that
$$\operatorname{Re}\left[-\langle x, y \rangle + \langle x, e \rangle \langle e, y \rangle\right] \leq \frac{1}{16} \cdot \frac{1}{\lambda (1 - \lambda)} |\Gamma - \gamma|^2$$

giving

(3.14) $$\operatorname{Re}\left[\langle x, y \rangle - \langle x, e \rangle \langle e, y \rangle\right] \geq -\frac{1}{16} \cdot \frac{1}{\lambda (1 - \lambda)} |\Gamma - \gamma|^2.$$

Consequently, by (3.3) and (3.14) we deduce the desired inequality (3.13). ∎

**Remark 2.** If $M, m \in \mathbb{R}$ with $M > m$ and, for $\lambda \in (0, 1)$,

$$\left\| \lambda x + (1 - \lambda) y - \frac{M + m}{2} e \right\| \leq \frac{1}{2} (M - m) \tag{3.15}$$

then

$$\langle x, y \rangle - \langle x, e \rangle \langle e, y \rangle \leq \frac{1}{16} \cdot \frac{1}{\lambda (1 - \lambda)} (M - m)^2.$$

If (3.15) holds with "$\pm$" instead of "$+$", then

$$|\langle x, y \rangle - \langle x, e \rangle \langle e, y \rangle| \leq \frac{1}{16} \cdot \frac{1}{\lambda (1 - \lambda)} (M - m)^2. \tag{3.16}$$

**Remark 3.** If $\lambda = \frac{1}{2}$ in (3.1) or (3.2), then we obtain the result from [2], i.e.,

$$\operatorname{Re} \left\langle \Gamma e - \frac{x + y}{2}, \frac{x + y}{2} - \gamma e \right\rangle \geq 0 \tag{3.17}$$

or, equivalently

$$\left\| \frac{x + y}{2} - \frac{\gamma + \Gamma}{2} e \right\| \leq \frac{1}{2} |\Gamma - \gamma| \tag{3.18}$$

implies

$$\operatorname{Re} [\langle x, y \rangle - \langle x, e \rangle \langle e, y \rangle] \leq \frac{1}{4} |\Gamma - \gamma|^2. \tag{3.19}$$

The constant $\frac{1}{4}$ is best possible in (3.19).

For $\lambda = \frac{1}{2}$, Corollary 4 and Remark 2 will produce the corresponding results obtained in [2]. We omit the details.

## 4. Integral Inequalities

Let $(\Omega, \Sigma, \mu)$ be a measure space consisting of a set $\Omega$, a $\sigma$-algebra of parts $\Sigma$ and a countably additive and positive measure $\mu$ on $\Sigma$ with values in $\mathbb{R} \cup \{\infty\}$. Denote by $L^2(\Omega, \mathbb{K})$ the Hilbert space of all real or complex valued functions $f$ defined on $\Omega$ and 2-integrable on $\Omega$, i.e.,

$$\int_\Omega |f(s)|^2 \, d\mu(s) < \infty.$$

The following proposition holds

**Proposition 1.** If $f, g, h \in L^2(\Omega, \mathbb{K})$ and $\varphi, \Phi, \gamma, \Gamma \in \mathbb{K}$, are so that $\operatorname{Re}(\Phi \overline{\varphi}) > 0$, $\operatorname{Re}(\Gamma \overline{\gamma}) > 0$, $\int_\Omega |h(s)|^2 \, d\mu(s) = 1$ and

$$\int_\Omega \operatorname{Re} \left[ (\Phi h(s) - f(s)) \left( \overline{f(s)} - \overline{\varphi} \overline{h(s)} \right) \right] d\mu(s) \geq 0 \tag{4.1}$$

$$\int_\Omega \operatorname{Re} \left[ (\Gamma h(s) - g(s)) \left( \overline{g(s)} - \overline{\gamma} \overline{h(s)} \right) \right] d\mu(s) \geq 0$$

or, equivalently

$$\left( \int_\Omega \left| f(s) - \frac{\Phi + \varphi}{2} h(s) \right|^2 d\mu(s) \right)^{\frac{1}{2}} \leq \frac{1}{2} |\Phi - \varphi|, \tag{4.2}$$

$$\left( \int_\Omega \left| g(s) - \frac{\Gamma + \gamma}{2} h(s) \right|^2 d\mu(s) \right)^{\frac{1}{2}} \leq \frac{1}{2} |\Gamma - \gamma|,$$

then we have the following Grüss type integral inequality

$$(4.3) \quad \left| \int_\Omega f(s) \overline{g(s)} d\mu(s) - \int_\Omega f(s) \overline{h(s)} d\mu(s) \int_\Omega h(s) \overline{g(s)} d\mu(s) \right|$$

$$\leq \frac{1}{4} M(\varphi, \Phi) M(\gamma, \Gamma) \left| \int_\Omega f(s) \overline{h(s)} d\mu(s) \int_\Omega h(s) \overline{g(s)} d\mu(s) \right|$$

where

$$M(\varphi, \Phi) := \left[ \frac{(|\Phi| - |\varphi|)^2 + 4\left[|\Phi \overline{\varphi}| - \operatorname{Re}(\Phi \overline{\varphi})\right]}{\operatorname{Re}(\Phi \overline{\varphi})} \right]^{\frac{1}{2}}.$$

The constant $\frac{1}{4}$ is best possible.

The proof follows by Theorem 3 on choosing $H = L^2(\Omega, \mathbb{K})$ with the inner product

$$\langle f, g \rangle := \int_\Omega f(s) \overline{g(s)} d\mu(s).$$

We omit the details.

**Remark 4.** *It is obvious that a sufficient condition for (4.1) to hold is*

$$\operatorname{Re}\left[ (\Phi h(s) - f(s)) \left( \overline{f(s)} - \overline{\varphi} \overline{h(s)} \right) \right] \geq 0,$$

*and*

$$\operatorname{Re}\left[ (\Gamma h(s) - g(s)) \left( \overline{g(s)} - \overline{\gamma} \overline{h(s)} \right) \right] \geq 0,$$

*for $\mu$-a.e. $s \in \Omega$, or equivalently,*

$$\left| f(s) - \frac{\Phi + \varphi}{2} h(s) \right| \leq \frac{1}{2} |\Phi - \varphi| |h(s)| \quad \text{and}$$

$$\left| g(s) - \frac{\Gamma + \gamma}{2} h(s) \right| \leq \frac{1}{2} |\Gamma - \gamma| |h(s)|,$$

*for $\mu$-a.e. $s \in \Omega$.*

The following result may be stated as well.

**Corollary 5.** *If $z, Z, t, T \in \mathbb{K}$, $\mu(\Omega) < \infty$ and $f, g \in L^2(\Omega, \mathbb{K})$ are such that:*

$$(4.4) \quad \operatorname{Re}\left[ (Z - f(s)) \left( \overline{f(s)} - \bar{z} \right) \right] \geq 0,$$

$$\operatorname{Re}\left[ (T - g(s)) \left( \overline{g(s)} - \bar{t} \right) \right] \geq 0 \quad \text{for a.e. } s \in \Omega$$

*or, equivalently*

$$(4.5) \quad \left| f(s) - \frac{z + Z}{2} \right| \leq \frac{1}{2} |Z - z|,$$

$$\left| g(s) - \frac{t + T}{2} \right| \leq \frac{1}{2} |T - t| \quad \text{for a.e. } s \in \Omega;$$

*then we have the inequality*

$$(4.6) \quad \left| \frac{1}{\mu(\Omega)} \int_\Omega f(s) \overline{g(s)} d\mu(s) - \frac{1}{\mu(\Omega)} \int_\Omega f(s) d\mu(s) \cdot \frac{1}{\mu(\Omega)} \int_\Omega \overline{g(s)} d\mu(s) \right|$$

$$\leq \frac{1}{4} M(z, Z) M(t, T) \left| \frac{1}{\mu(\Omega)} \int_\Omega f(s) d\mu(s) \cdot \frac{1}{\mu(\Omega)} \int_\Omega \overline{g(s)} d\mu(s) \right|.$$

**Remark 5.** *The case of real functions incorporates the following interesting inequality*

$$
(4.7) \quad \left| \frac{\mu(\Omega) \int_\Omega f(s) g(s) d\mu(s)}{\int_\Omega f(s) d\mu(s) \int_\Omega g(s) d\mu(s)} - 1 \right| \leq \frac{1}{4} \cdot \frac{(Z-z)(T-t)}{\sqrt{ztZT}}
$$

*provided $\mu(\Omega) < \infty$,*

$$z \leq f(s) \leq Z, t \leq g(s) \leq T$$

*for $\mu-a.e.\ s \in \Omega$, where $z, t, Z, T$ are real numbers and the integrals at the denominator are not zero. Here the constant $\frac{1}{4}$ is best possible in the sense mentioned above.*

Using Theorem 3 we may state the following result as well.

**Proposition 2.** *If $f, g, h \in L^2(\Omega, \mathbb{K})$ and $\gamma, \Gamma \in \mathbb{K}$ are such that $\int_\Omega |h(s)|^2 d\mu(s) = 1$ and*

$$
(4.8) \quad \int_\Omega \{\operatorname{Re}[\Gamma h(s) - (\lambda f(s) + (1-\lambda) g(s))]
$$
$$
\times \left[\lambda \overline{f(s)} + (1-\lambda) \overline{g(s)} - \bar{\gamma} \bar{h}(s)\right]\} d\mu(s)
$$
$$
\geq 0
$$

*or, equivalently,*

$$
(4.9) \quad \left( \int_\Omega \left| \lambda f(s) + (1-\lambda) g(s) - \frac{\gamma + \Gamma}{2} h(s) \right|^2 d\mu(s) \right)^{\frac{1}{2}} \leq \frac{1}{2} |\Gamma - \gamma|,
$$

*then we have the inequality*

$$
(4.10) \quad I := \int_\Omega \operatorname{Re}\left[f(s) \overline{g(s)}\right] d\mu(s)
$$
$$
- \operatorname{Re}\left[ \int_\Omega f(s) \overline{h(s)} d\mu(s) \cdot \int_\Omega h(s) \overline{g(s)} d\mu(s) \right]
$$
$$
\leq \frac{1}{16} \cdot \frac{1}{\lambda(1-\lambda)} |\Gamma - \gamma|^2.
$$

*The constant $\frac{1}{16}$ is best possible.*

*If (4.8) and (4.9) hold with "$\pm$" instead of "$+$" (see Corollary 4), then*

$$
(4.11) \quad |I| \leq \frac{1}{16} \cdot \frac{1}{\lambda(1-\lambda)} |\Gamma - \gamma|^2.
$$

**Remark 6.** *It is obvious that a sufficient condition for (4.8) to hold is*
(4.12)
$$
\operatorname{Re}\left\{[\Gamma h(s) - (\lambda f(s) + (1-\lambda) g(s))] \cdot \left[\lambda \overline{f(s)} + (1-\lambda) \overline{g(s)} - \bar{\gamma} \bar{h}(s)\right]\right\} \geq 0
$$

*for a.e. $s \in \Omega$, or equivalently*

$$
(4.13) \quad \left| \lambda f(s) + (1-\lambda) g(s) - \frac{\gamma + \Gamma}{2} h(s) \right| \leq \frac{1}{2} |\Gamma - \gamma| |h(s)|
$$

*for a.e. $s \in \Omega$.*

Finally, the following corollary holds.

**Corollary 6.** *If $Z, z \in \mathbb{K}$, $\mu(\Omega) < \infty$ and $f, g \in L^2(\Omega, \mathbb{K})$ are such that*

(4.14) $\quad \operatorname{Re}\left[\left(Z - (\lambda f(s) + (1-\lambda)g(s))\right)\left(\lambda \overline{f(s)} + (1-\lambda)\overline{g(s)} - \overline{z}\right)\right] \geq 0$

*for a.e. $s \in \Omega$, or, equivalently*

(4.15) $\quad \left|\lambda f(s) + (1-\lambda)g(s) - \dfrac{z+Z}{2}\right| \leq \dfrac{1}{2}|Z - z|,$

*for a.e. $s \in \Omega$, then we have the inequality*

$$J := \dfrac{1}{\mu(\Omega)}\int_\Omega \operatorname{Re}\left[f(s)\overline{g(s)}\right]d\mu(s)$$
$$- \operatorname{Re}\left[\dfrac{1}{\mu(\Omega)}\int_\Omega f(s)\,d\mu(s) \cdot \dfrac{1}{\mu(\Omega)}\int_\Omega \overline{g(s)}d\mu(s)\right]$$
$$\leq \dfrac{1}{16} \cdot \dfrac{1}{\lambda(1-\lambda)}|Z-z|^2.$$

*If (4.14) and (4.15) hold with " $\pm$ " instead of " $+$ ", then*

(4.16) $\quad |J| \leq \dfrac{1}{16} \cdot \dfrac{1}{\lambda(1-\lambda)}|Z-z|^2.$

**Remark 7.** *It is obvious that if one chooses the discrete measure above, then all the inequalities in this section may be written for sequences of real or complex numbers. We omit the details.*

## References

[1] S.S. Dragomir, A generalization of Grüss' inequality in inner product spaces and applications, *J. Math. Anal. Appl.*, **237**(1999), 74-82.

[2] S.S. Dragomir, Some Grüss' type inequalities in inner product spaces, accepted in *J. Ineq. Pure & Appl. Math.* (`http://jipam.vu.edu.au`).

[3] S.S. Dragomir and I. Gomm, Some integral and discrete versions of the Grüss inequality for real and complex functions and sequences, *RGMIA Res. Rep. Coll.*, **5**(2003), No. 3, Article 9 [ON LINE `http://rgmia.vu.edu.au/v5n3.html`]

School of Computer Science and Mathematics, Victoria University of Technology, PO Box 14428, Melbourne City MC 8001, Victoria, Australia.
*E-mail address*: `sever.dragomir@vu.edu.au`
*URL*: `http://rgmia.vu.edu.au/SSDragomirWeb.html`